\documentclass[11pt]{amsart}

\usepackage[left=3cm,top=3cm,right=3cm,bottom=3cm]{geometry}
\usepackage{amssymb, amsmath, amsthm, units}
\usepackage{mathtools}
\usepackage{verbatim}
\usepackage{comment}
\usepackage{graphicx}
\usepackage{enumerate}
\usepackage{color}
\usepackage{hyperref}
\usepackage[capitalize,nameinlink,noabbrev]{cleveref}
\usepackage{float}

\DeclareGraphicsExtensions{.pdf,.jpg,.png}

\theoremstyle{plain}
\newtheorem{theorem}{Theorem}

\newtheorem{proposition}[theorem]{Proposition}

\newtheorem{conjecture}[theorem]{Conjecture}

\theoremstyle{definition}
\newtheorem{definition}[theorem]{Definition}

\newcommand{\eps}{\varepsilon}

\DeclarePairedDelimiter\abs{\lvert}{\rvert}
\makeatletter
\let\oldabs\abs
\def\abs{\@ifstar{\oldabs}{\oldabs*}}

\begin{document}

\title{Intersecting hexagons in 3-space}

\author{J\'{o}zsef Solymosi}
\address{\noindent Department of Mathematics, University of British Columbia, Vancouver, BC, Canada V6T 1Z2}
\email{solymosi@math.ubc.ca}
\thanks{The first author was supported by NSERC, NKFI KKP 133819, and OTKA NK 104183 grants}

\author{Ching Wong}
\address{\noindent Department of Mathematics, University of British Columbia, Vancouver, BC, Canada V6T 1Z2} \email{ching@math.ubc.ca}

\date{}

\begin{abstract}
Two hexagons in the space are said to intersect {\em heavily} if their intersection consists of at least one common vertex as well as an interior point. We show that the number of hexagons on $n$ points in 3-space without heavy intersections is $o(n^2)$, under the assumption that the hexagons are `fat'.
\end{abstract}

\maketitle

\section{Introduction}

The problem of finding the maximum number of hyperedges in a geometric hypergraph in $d$-dimensional space with certain forbidden configurations (intersections) is a general problem in discrete geometry. Several such problems were considered by Dey and Pach in \cite{DeyPac98}. We are interested in finding the maximum number of (planar and convex) polygons on some vertex set of $n$ points in $3$-space, where no two of them are allowed to intersect in certain ways. In this paper, we confine our attention to 3-space, and by polygons we mean planar polygons, i.e. the vertices are co-planar, which are convex. As usual, a $k$-$gon$ (where $k \geq 3$) is a polygon with $k$ vertices.

\subsection{Almost disjoint polygons}

\begin{figure}[ht!]
	\centering
	\includegraphics[width=75mm]{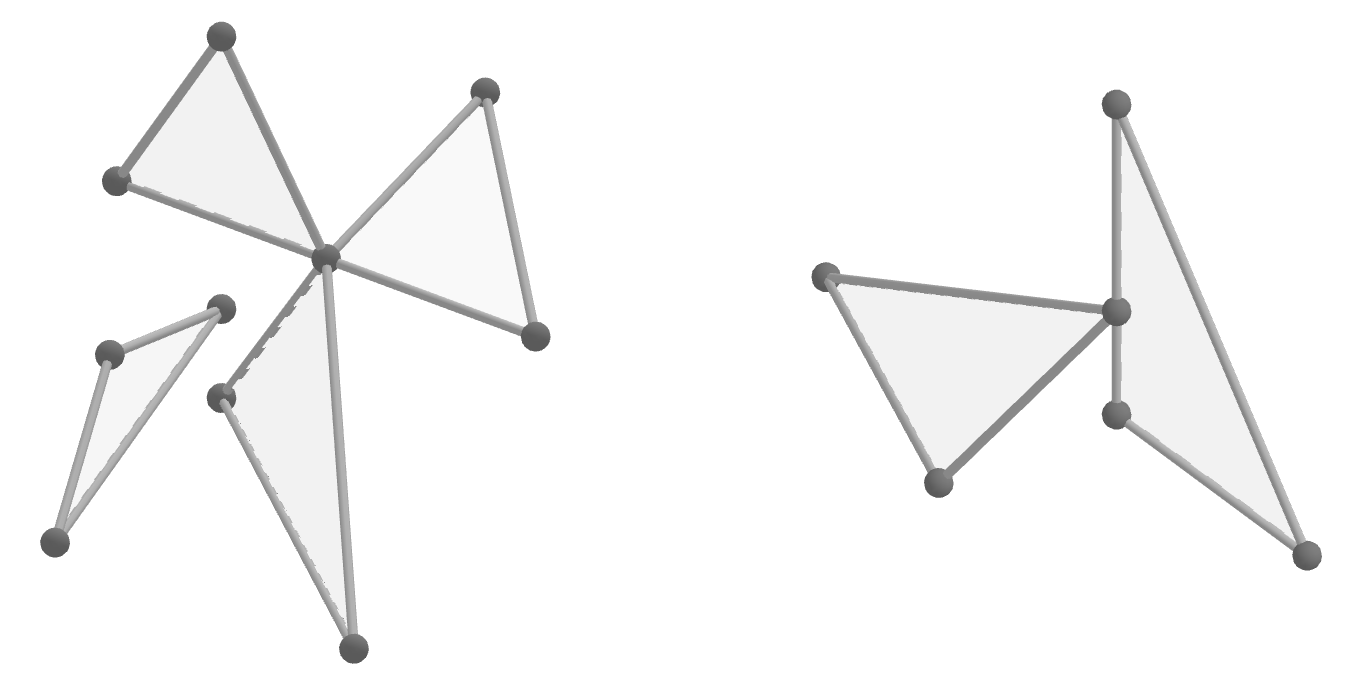}
	\caption{(Left) A set of pairwise almost disjoint triangles. (Right) Two triangles that are not almost disjoint.}
	\label{fig:almostdisjoint}
\end{figure}

It was asked by Gil Kalai \cite{Kal} and independently by G\"unter Ziegler (quoted from \cite{KarSol02}) what the maximum possible number of triangles spanned by $n$ points is, such that any two are almost disjoint:

\begin{definition}[Almost disjoint polygons]
	Two polygons in $3$-space are said to be {\em almost disjoint} if they are either disjoint, or their intersection consists of one common vertex. See \Cref{fig:almostdisjoint}.
\end{definition}

Let $f_1(n,k)$ be the maximum possible number of pairwise almost disjoint $k$-gons on $n$ points in $3$-space. It is easy to see that $f_1(n,k) \geq f_1(n,k+1)$ for all $k \geq 3$, by arbitrarily forming a $k$-gon from each $(k+1)$-gon.

We remark that a set of pairwise almost disjoint triangles on $n$ points gives a partial Steiner triple system of order $n$. A simple counting reveals that such system has at most $n(n-1)/6$ triples. It follows that
\[
f_1(n,k) \leq f_1(n,3) \leq \dfrac{n(n-1)}{6},
\]
for every $k \geq 3$. 

K\'{a}rolyi and Solymosi \cite{KarSol02} constructed configurations showing that $f_1(n,3) \geq cn^{3/2}$ for some universal constant $c>0$. Finding sharper lower bounds seems like a very hard problem. In fact, it is not even known if the genus of a polytope on $n$ vertices can have order ${n^2}$. If so, the magnitude of $f_1(n,3)$ would be $n^2$. The best lower bound of the largest genus is $n \log n$, due to a construction of McMullen, Schulz and Wills \cite{McmSchWil83}. For more details, we refer the interested reader to \cite{Zie08} where Ziegler gives a simplified construction providing the same bound.

When $k=3$, a tight asymptotic bound can be obtained if we relax the assumption a little by allowing the triangles to intersect also in an (entire) edge. Let $f_2(n,k)$ be the maximum possible number of $k$-gons on $n$ points so that any two of them are either disjoint, or intersect in a vertex or an edge. Since we allow more intersections, we have $f_2(n,k) \geq f_1(n,k)$ for all $n$ and $k \geq 3$. It is again clear that $f_2(n,k)$ decreases with $k$, for all $n$.

\begin{proposition}
	\label{prop:f2}
	As $n\to\infty$,
	\[
	f_2(n,3)=\Theta(n^2).
	\]
\end{proposition}

The proof is given in \Cref{sec:f2proof}.

\subsection{Non-heavily intersecting polygons}

In this chapter, we focus on an even more relaxed assumption on the intersections of the $k$-gons.

\begin{definition}[heavily intersecting polygons]
	Two polygons in $3$-space are said to {\em intersect heavily} if their intersection consists of at least one common vertex as well as an interior point.
\end{definition}

\begin{figure}[H]
	\centering
	\includegraphics[width=75mm]{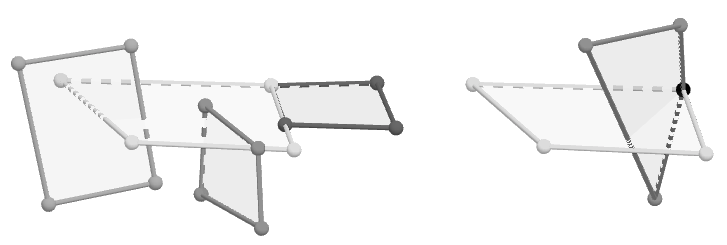}
	\caption{(Left) A set of 4 quadrilaterals without heavy intersections. (Right) Two quadrilaterals that intersect heavily.}
	\label{fig:heavilyIntersect}
\end{figure}

We say that a collection of $k$-gons has no heavy intersections if no two of these polygons intersect heavily. See \Cref{fig:heavilyIntersect}.

Let $f_3(n,k)$ be the maximum possible number of $k$-gons without heavy intersections on $n$ points in $3$-space. It is, once again, true that $f_3(n,k) \geq f_3(n,k+1)$ and $f_3(n,k) \geq f_2(n,k)$, for all $k \geq 3$.

In such arrangements, two $k$-gons cannot share a diagonal and so $f_3(n,k) = O(n^2)$ for $k \geq 4$. In fact, the proof of \Cref{prop:f2} (first part) works here as well, and so
\[
f_3(n,k) < n^2,
\]
for every $k \geq 3$.

\begin{figure}[H]
	\centering
	\includegraphics[width=50mm]{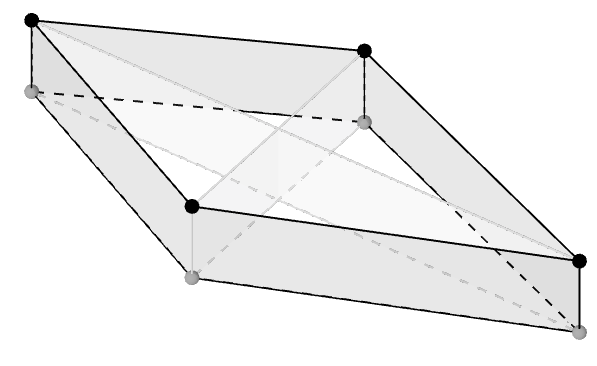}
	\caption{A set of $\binom{8/2}{2} = 6$ quadrilaterals without heavy intersections.}
	\label{fig:f3LowerBound}
\end{figure}

This upper bound is actually sharp, in magnitude, for triangles and quadrilaterals ($k=3,4$). One can give a construction of $\Omega(n^2)$ quadrilaterals on $n$ points without heavy intersections as follows: Let $n$ be an even number and suppose we are given $n/2$ points $P_1, \dots, P_{n/2}$ in general position (no three points collinear) on a plane $\pi$. Fix any vector $v$ not parallel to $\pi$. Then the $n$ points $P_1, \dots, P_{n/2}, P_1+v, \dots, P_{n/2}+v$ are incident to $\binom{n/2}{2} = (n^2-2n)/8$ desired quadrilaterals with vertices $P_i, P_j, P_j+v, P_i+v$, where $1 \leq i < j \leq n/2$. \Cref{fig:f3LowerBound} shows an example when $n=8$.

When $k=6$, we show that the number of hexagons without heavy intersections is $o(n^2)$, under an extra assumption on the `fatness' of the hexagons defined below.

\begin{definition}[Fat hexagons]\label{defn:fat}
	Let $c\geq1$ and $0<\alpha<\pi/2$. A hexagon is $(c,\alpha)$-{\em fat} if
	\begin{enumerate}
		\item the ratio of any two sides is bounded between $1/c$ and $c$, and
		
		\item it has three non-neighbour vertices having interior angles between $\alpha$ and $\pi-\alpha$.
	\end{enumerate}
\end{definition}

Our main tool is the Triangle Removal Lemma of Ruzsa and Szemer\'edi, which states that for any $\eps>0$, there exists $\delta>0$ such that any graph on $n$ vertices with at least $\eps{n^2}$ pairwise edge-disjoint triangles, has at least $\delta{n^3}$ triangles in total. See \cite{RuzSze78} for the original formulation of this result. The following is the precise statement of our main theorem, which is proved in \Cref{sec:hexagonProof}.

\begin{theorem}
	\label{thm:hexagon}
	For any $c\geq1$ and $0<\alpha<\pi/2$, there is a function $F_{(c,\alpha)}(n),$
	\[\frac{F_{(c,\alpha)}(n)}{n^2}\rightarrow 0 \quad \text{    as    } \quad n\rightarrow \infty,\]
	such that any family of $(c,\alpha)$-fat hexagons in $3$-space on $n$ points without heavy intersections has size at most  $F_{(c,\alpha)}(n)$.
\end{theorem}


\section{Proof of \Cref{prop:f2}}
\label{sec:f2proof}

\begin{figure}[H]
	\centering
	\includegraphics[width=75mm]{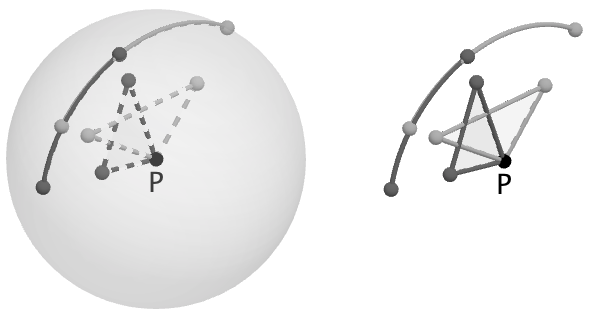}
	\caption{The graph $G$ does not have multiple edges, or there are some triangles intersecting in an unwanted way.}
	\label{fig:projectionSphere1}
\end{figure}

\begin{figure}[ht]
	\centering
	\includegraphics[width=75mm]{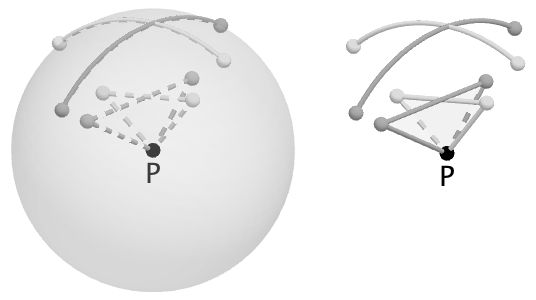}
	\caption{If there were two triangles both incident to $P$ so that the projections of their corresponding line segments intersect internally, then these two triangles would intersect in an unwanted way.}
	\label{fig:projectionSphere2}
\end{figure}

We first show that $f_2(n,3) < n^2$. Given such a set of triangles on $n$ points. Pick any of these $n$ points, say $P$, and project the remaining $n-1$ points in $3$-space onto a sphere $S$ centred at $P$. We want to upper bound the number of triangles incident with the point $P$. For each of such triangles, say $PQR$, we project the line segment between $Q$ and $R$ onto the sphere $S$. These geodesic segments, together with the projected points, form a graph $G$ on $S$ having at most $n-1$ vertices. Here, we subdivide an edge (geodesic segment) if there are vertices lying on it. We note that if there were multiple edges on $G$, then their corresponding triangles would lie on the same plane and intersect in an interior point, as shown in \Cref{fig:projectionSphere1}. Hence, the number of edges in $G$ is at least the number of triangles incident to $P$. As illustrated in \Cref{fig:projectionSphere2}, the graph $G$ is planar, and so it has at most $3(n-1)-6 = 3n-9$ edges. Hence,
\[
f_2(n,3) \leq \dfrac{n}{3} (3n-9) < n^2.
\]

\begin{figure}[ht]
	\centering
	\includegraphics[width=75mm]{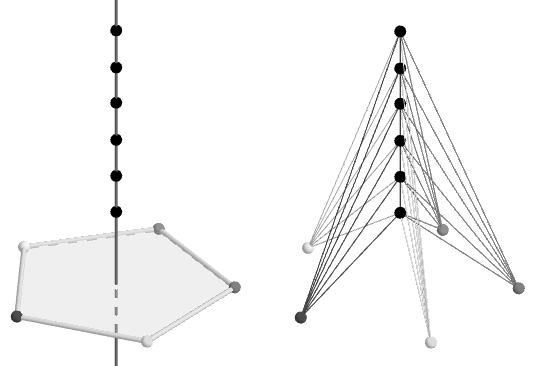}
	\caption{(Left) There are $m=5$ points on the $xy$-plane and $m+1=6$ points on the $z$-axis. (Right) There are $m^2=25$ triangles on $2m+1=11$ points.}
	\label{fig:christmasTree}
\end{figure}

Now, we show that $f_2(n,3) \geq (n-1)^2/4$ whenever $n$ is odd. To see this, we are using the well-known Christmas tree arrangement. Let there be $m$ points on a circle centred at the origin on the $xy$-plane  and $m+1$ points on $z$-axis, as in \Cref{fig:christmasTree} (left). So we have a total of $n=2m+1$ points in $3$-space. We consider the $m^2=(n-1)^2/4$ triangles with one vertex on the circle, and the other two vertices being a consecutive pair of points chosen on the $z$-axis. See \Cref{fig:christmasTree} (right). It is easy to see that if two triangles are not disjoint, they intersect in either a vertex or an edge.

\section{Proof of \Cref{thm:hexagon}}
\label{sec:hexagonProof}
Let $\eps>0$. Given $\eps{n^2}$ many $(c,\alpha)$-fat hexagons on $n$ vertices in $3$-space. We show that two of these hexagons intersect heavily, when $n$ is large enough. We may assume, in particular, that no two hexagons share a diagonal.

To reduce the dimension of the ambient space, we project these hexagons onto a random plane such that no two vertices share the same projection and that a positive fraction of the hexagons is $(c',\alpha')$-fat. Indeed, if we project a $(c,\alpha)$-fat hexagon $H$ onto a plane making an angle at most $\theta<\pi/2$ with the plane containing $H$, it is straightforward to show that the projected hexagon is $(c',\alpha')$-fat, where
$$c'=\dfrac{c}{\cos\theta}\quad\quad\text{and}\quad\quad\alpha'=\cos^{-1}\left(\dfrac{\cos\alpha+\sin^2\theta}{\cos^2\theta}\right).$$

The existence of heavily intersecting hexagons relies on a similar-slope property. This can be described quantitatively by the difference of  two angles of inclination. To this end, let $\phi>0$ be the smallness of such differences which is to be determined later.

We choose from the $(c',\alpha')$-fat projected hexagons the most popular family consisting of $\eps'n^2$ hexagons, which have inscribed triangles of similar shapes and orientations.

More precisely, let us enumerate by any order the projected hexagons as $\{H_i\}$ and label their vertices as $A_i,B_i,C_i,D_i,E_i,F_i$, oriented counter-clockwise, where $B_i,D_i,F_i$ are the three non-neighbour vertices having angles between $\alpha'$ and $\pi-\alpha'$.

There exists a positive fraction of these hexagons so that for any $i,j$, the inclined angles of the diagonals $A_iC_i$ and $A_jC_j$ differ by at most $\phi$. Similarly the same property holds true for the diagonals $C_iE_i$ and $E_iA_i$ in yet a sub-collection of $\eps'n^2$ hexagons.

We define $G$ to be the graph whose vertices are the $n$ projected points and whose edges are from the triangles formed by the vertices $A_i,C_i,E_i$ chosen above. Then, as we assumed that no two hexagons share a diagonal, $G$ contains $\eps'n^2$ edge-disjoint triangles. An application of the Triangle Removal Lemma yields, when $n$ is large enough, a triangle $T$ whose edges come from three different hexagons, say $H_1$, $H_2$ and $H_3$. For each $i=1,2,3$, let $T_i$ be the triangle $A_iC_iE_i$.

We are ready to study the intersection properties of these three hexagons in the $3$-space. In other words, we now `unproject' the $n$ points.

\begin{figure}[H]
	\centering
	\includegraphics[width=60mm]{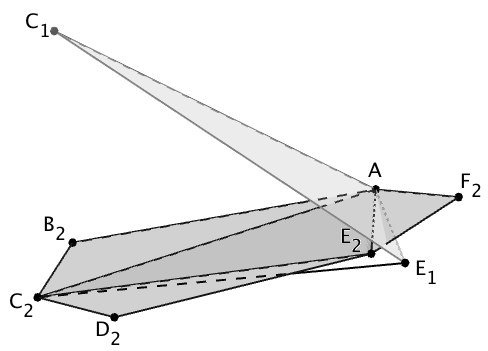}
	\caption{The triangle $T_1=AC_1E_1$ and the hexagon $H_2=AB_2C_2D_2E_2F_2$ intersect heavily. Here the triangle $T$ is $AC_2E_1$.}
	\label{heavy}
\end{figure}

Two of the triangles, say $T_1$ and $T_2$, lie on the same side of $T$ and let $T_1$ be the triangle making a larger angle with $T$. Then, as shown in \Cref{heavy}, the hexagon $H_2$ intersects heavily with the triangle $T_1$, and hence with the hexagon $H_1$, as long as the three non-neighbour vertices $B_1,D_1,F_1$ lie outside of the triangle $T$ on the plane of projection, which is guaranteed if we choose
$$\phi<\tan^{-1}\left(\dfrac{\sin\alpha'}{c'+\cos\alpha'}\right),$$
the right hand side being a lower bound of the six angles $B_1A_1C_1$ etc. under the $(c',\alpha')$-fatness assumption. This completes the proof.

\section{Open problems}
The original question of Kalai and Ziegler is still widely open even if we restrict the question to fat triangles. Let us say that a triangle is fat if all its angles are between 50 and 70 degrees. Let us state a special case of the Kalai-Ziegler problem
\begin{conjecture}
In 3-space the number of almost disjoint fat triangles spanned by $n$
 points is $o(n^2).$
\end{conjecture}

One would expect a $O(n^{2-c})$ type upper bound with some $c>0,$ but we can't even show $o(n^2).$ On the other hand we are not aware of any construction with a superlinear number of almost disjoint fat triangles. 

Although Regularity Lemmas were used in discrete geometry, it is very likely that in the proof of Theorem  \ref{thm:hexagon} one could substitute it with some geometric arguments providing much better bounds.

\end{document}